\documentclass[journal,9pt]{IEEEtran}
\usepackage{algorithmic,url,color}
\usepackage{algorithm}
\usepackage{graphicx}
\usepackage{cite}
\usepackage{amssymb}
\usepackage{amsmath,amsfonts,mathrsfs}
\usepackage{mathabx}
\usepackage{textcomp}

\newtheorem{Lemma}{Lemma}

\newtheorem{Definition}{Definition}
\newtheorem{Remark}{Remark}

\newtheorem{Theorem}{Theorem}

\newcommand{\Rank}{\mathrm{Rank}}

\newcommand{\V}{\mathcal{V}}
\newcommand{\IM}{\mathrm{im}}
\newcommand{\LMAP}{\mathcal{M}}
\newcommand{\LOPT}{\mathcal{M}_{opt}}

\newcommand{\obs}[1]{\mathcal{O}_U(#1,y)}
\newcommand{\mobs}[1]{\mathscr{O}_{\widecheck U}(y)(#1)}

\newcommand{\Id}{\mathcal{I}}
\newcommand{\diag}{\mathrm{diag}}
\newcommand{\rank}{\mathrm{rank}}

\begin{document}
\title{Minimax Observers for Linear DAEs}
\date{}
\author{Sergiy~Zhuk and  Mihaly~Petreczky
\thanks{S. Zhuk is with IBM Research - Ireland, Damastown, Dublin 15, Ireland, \texttt{sergiy.zhuk@ie.ibm.com}}
\thanks{CNRS, Centrale Lille, UMR 9189 - CRIStAL-Centre de Recherche en Informatique, Signal et Automatique de Lille, F-59000 Lille, France,  \texttt{mihaly.petreczky@ec-lille.fr}}
}
\maketitle
\begin{abstract}
In this note we construct minimax observers for linear stationary DAEs with bounded uncertain inputs, given noisy measurements.
We prove a new duality principle and show that a finite (infinite) horizon minimax observer exists if and only if the DAE is $\ell$-impulse observable ($\ell$-detectable) . Remarkably, the regularity of the DAE is not required.
\end{abstract}

\vspace{-10pt}
\section{Introduction}
\label{sec:introduction}
Consider a linear Differential-Algebraic Equation (DAE): 
\vspace{-5pt}
 \begin{equation}
\label{eq:dae}
 \dfrac {d(Fx(t))}{dt}=Ax(t)+f(t)\,, y(t) = Hx(t) +\eta(t)\,, 
\vspace{-5pt}
 \end{equation}
where $F, A\in\mathbb R^{m \times n}$ and $H\in\mathbb R^{p\times n}$ are given matrices and $f$, $\eta$ are deterministic noises. We will be interested in solutions for which $(Fx(0),f,\eta)$ are \textit{a priori} unknown elements of a bounding set $\mathscr E_\infty:=\bigcap_{t_1>0}\{(x_0,f,\eta): \rho(x_0,f,\eta,t_1,Q_0,Q,R)\le 1\}$ where
\vspace{-5pt}
\begin{equation}
 \label{eq:rhox0feta0}
 \begin{split}
 & \rho(x_0,f,\eta,t_1,Q_0,Q,R):=x_0^{\top} Q_0 x_0 +  \\
  & + \int_{0}^{t_1} (f^{\top}(t) Qf(t) + \eta^{\top}(t) R\eta(t)) dt\,,
 \end{split}
\end{equation}
and $Q_0,Q,R$ are symmetric positive definite matrices. Intuitively, $(Fx(0),f,\eta)\in \mathscr E_\infty$ implies that the initial state $Fx(0)$ is bounded and the noise signals $(f,\eta)$ have bounded energy, i.e. bounded $L^2$ norm. The latter is a standard assumption in robust control~\cite{RobustControl}.

In this paper we extend the classical results on minimax observers (see~\cite{KrenerMinimax,Tempo1985,Nakonechnii1978,RobustControl}) to DAEs. Namely, we estimate a linear function $t\mapsto Lx(t)$, $L \in \mathbb{R}^{1 \times n}$ of the state trajectory of \eqref{eq:dae}, based on the output $y$. Since~\eqref{eq:dae} contain uncertain terms $f,\eta$, it follows that, even if $F$ was invertible, there could be various state trajectories $x$ of \eqref{eq:dae}, and hence several values of $Lx(t)$, which are consistent with $y$. For the case of rectangular $F$, the DAE~\eqref{eq:dae} may have several solutions. This introduces an additional source of uncertainty. Hence, the best we can do 
is to construct a function $t\mapsto \widehat{Lx}(t)$ such that \textbf{for all} state trajectories $x$ of \eqref{eq:dae}, which are consistent with $y$ and correspond to $(Fx(0),f,\eta)\in \mathscr E_\infty$, $Lx(t)$ lies in the interval $[\widehat{Lx}(t)-r,\widehat{Lx}(t)+r]$,
and $r$ is as small as possible. A dynamical system whose output is $\widehat{Lx}$ and whose input is $y$ is called a \emph{minimax observer}.
%
In what follows we take $Lx\!\!=\!\!\ell^{\top}Fx$ for some $\ell \in \mathbb{R}^m$. The reason for this is that the state vector $x$ of \eqref{eq:dae} may contain components which are free, possibly discontinuous, exogenous variables.
The estimation of the latter is problematic even for LTI systems \cite{Hautus1983}, if the observer is required to be an LTI system. 
In contrast, $Fx$ is always continuous and it is subject to a differential equation.
\subsubsection{Contribution}
Our main contribution is a new duality principle (Theorem~\ref{p:5}). This result allows us (i) to find necessary and sufficient conditions for existence of a minimax observer $\widehat{Lx}$ for finite and infinite time horizons (Theorems~\ref{p:1:theo},\ref{p:5}), and (ii) to show that $\widehat{Lx}$ can be represented as an output of a time-varying linear system in the finite horizon case (Theorem~\ref{t:observer-fin-hor}), and as an output of a stable LTI system in the case of infinite time horizon (Theorem~\ref{t:observer-infin-hor}).
Finally we prove that for bounded $f$ and $\eta$ with possibly unbounded $L^2$ norm, the estimation error remains bounded. 
\subsubsection{Motivation}
\label{sec:motivation}
The need for estimating the state of a general DAE arises in many applications, as it has been demonstrated by various authors, e.g., \cite{BergerThesis,DAEBookChapter,Darouach2009,Xu2007,Zhang2013,Darouach2013,HouMuelle1999,IshiharaTerra}. In this paper we extend the classical minimax framework to DAEs of the form~\eqref{eq:dae}. One motivation 
is to provide control-theoretic grounds for solving an \emph{issue of approximation and state estimation for Partial Differential Equations} (PDEs), described in~\cite{ZhukSISC13}. This issue may be resolved by using our results to incorporate PDE's approximation error into the estimation process, see  Section~\ref{sec:numerical-example} for an example.

\subsubsection{Related work}
\label{sec:related-work}
To the best of our knowledge, the results of this paper are new. Its preliminary version appeared in \cite{ZhukPetreczkyCDC13}. With respect to \cite{ZhukPetreczkyCDC13} the main differences are: (i) new (necessary and sufficient) conditions for existence of minimax observers (Theorem~\ref{p:5}), (ii) detailed proofs.
Duality principle for non-stationary DAEs was introduced in~\cite{Zhuk2012AMO}, provided $Fx(0)=0$. It was then used to derive a sub-optimal observer. In contrast, our duality theorems hold for DAEs with uncertain $Fx(0)=x_0$ and the constructed observers are optimal in that the worst-case estimation errors associated with the observers are minimal (see Definitions~\ref{d:finhorobs}-\ref{d:infhorobs}). The algorithm of~\cite{Zhuk2012sysid} constructing a finite horizon minimax observer by using projectors $FF^+$ is a special case of the one of this paper. Minimax observers for discrete time DAEs were considered in~\cite{Zhuk2009c}.

We note that some papers (see~\cite{DuanTut} for an overview) derive observers for regular DAEs by converting it to the Weierstrass canonical form. If DAE has index $j$, then the resulting ODE depends on the derivatives of the noise of order $j-1$, and so the observer would require bounds for them~\cite{Gerdin2007}.
In contrast, 
our observer is independent of DAE's index and works for bounded $f$ and $\eta$. The papers~\cite{Xu2007,Darouach2009} present design of $H_{\infty}$-observers which force $H_{\infty}$-norm of the ``error system'' to be less than a given number. In contrast, (i) minimax observers of this paper minimize a different error measure, (ii) unlike~\cite{Xu2007} we allow for non-regular DAEs, and (iii) unlike~\cite{Darouach2009} we present necessary and sufficient conditions for observer existence ($\ell$-detectability, Theorem~\ref{p:5}). Sufficient existence conditions (partial impulse observability) for the asymptotic functional observer for~\eqref{eq:dae} were proposed in~\cite{Darouach2013} provided $f=0,\eta=0$. In contrast, we deal with noisy systems and present necessary and sufficient conditions. Design of unknown input observers for DAEs was considered in~\cite{Zhang2013}. Since we deal with bounded unknown inputs and estimate just a part of DAE's state, our results do not compare directly with~\cite{Zhang2013}. In this paper we consider deterministic noise, stochastic or fuzzy uncertainties have been addressed for example in \cite{FuzzyDAE,Pulch}.

\subsubsection{Outline of the paper}
\label{sec:outline}
The definitions of minimax observers are formulated in Section~\ref{sec:problem-statement}. The main results are presented in
Section~\ref{sec:main}. The application of the observers to PDEs is given in Section~\ref{sec:numerical-example}. The proofs are given in Section \ref{sect:proof}.

\subsubsection{Notation}
\label{sec:notation-1}
$\Id_n$ denotes the $n \times n$ identity matrix;
for a matrix $A$,  $A^+$ denotes its Moore-Penrose pseudoinverse.
Let  $I:=[0,t]$, $0 < t \in\mathbb{R}$ or $I:=[0,+\infty)$. Denote
by $AC(I,\mathbb{R}^n)$, $L^2(I,\mathbb{R}^{n})$, $L^2_{loc}(I,\mathbb{R}^n)$  respectively the set of all absolutely continuous, the set of all square integrable, and the set of all locally square integrable functions of the form $f:I \rightarrow \mathbb{R}^{n}$.
Recall that  $f:I \rightarrow \mathbb{R}^{n}$
is locally square integrable if its restriction to any compact interval $I^1 \subseteq I$
is square integrable.
We write $AC(I)$, $L^2(I)$ and $L^2_{loc}(I)$ referring to  $AC(I,\mathbb{R}^n)$, $L^2(I,\mathbb{R}^{n})$ and $L^2_{loc}(I,\mathbb{R}^{n})$ respectively, when $\mathbb{R}^n$ is clear from the context.
We say that $f(t)=g(t)$ holds almost everywhere (a.e.) on $I$ if $f(t)\ne g(t)$ only on a subset of $I$ of measure zero.
We identify linear time invariant (LTI) systems $\dot x=Ax+Bu$, $y=Cx+Du$, with
the tuple of matrices $(A,B,C,D)$.
By convention, the infimum of an empty set will be $+\infty$. Let $I=[0,t_1]$, $0 < t_1 \in \mathbb{R}$. For any $v \in L^2(I,\mathbb{R}^k)$, define  $\delta_{t_1}(v) \in L^2(I,\mathbb{R}^k)$ by $\delta_{t_1}(v)(s)=v(t_1-s)$ a.e.
For any  $\tau > 0$, positive definite matrices $S_0,S_1 \in \mathbb{R}^{m \times m}, S_2 \in \mathbb{R}^{p \times p}$, vector $g_0 \in \mathbb{R}^{m}$, and functions $g \in L^2_{loc}(I,\mathbb{R}^{m}), h \in L^2_{loc}(I,\mathbb{R}^{p})$, where $I$ contains $[0,\tau]$, define
\begin{equation} 
\vspace{-5pt}
 \label{eq:rhox0feta}
 \begin{split}
   & \rho(g_0,g,h,\tau,S_0,S_1,S_2):= g_0^{\top} S_0 g_0 + \\ & + \int_{0}^{\tau} (g^{\top}(t) S_1g(t) + h^{\top}(t) S_2 h(t)) dt.
 \end{split}
\end{equation}
The quantity $\rho(x_0,f,\eta,t_1,Q_0,Q,R)$ defined in \eqref{eq:rhox0feta0} is a special case of the notation defined
in \eqref{eq:rhox0feta} applied to $\tau=t_1$, $g_0=x_0$, $f=g$, $h=\eta$, $Q_0=S_0$, $Q=S_1$, $R=S_2$.

\section{Problem statement}
\label{sec:problem-statement}
We begin with a precise definition of DAE's solution on the interval $I$ of the form $I=[0,t_1]$, $0<t_1 \in \mathbb{R}$ or $I=[0,+\infty)$.
\begin{Definition}\label{sol:behav}
 A \emph{solution} of \eqref{eq:dae} on $I$ is a tuple $(x,f,y,\eta) \in L^2_{loc}(I,\mathbb{R}^n) \times L^2_{loc}(I,\mathbb{R}^m) \times L^2_{loc}(I,\mathbb{R}^p) \times L^2_{loc}(I,\mathbb{R}^p)$ such that: $Fx$ is absolutely continuous and $Fx(t)=Fx(0)+\int_{0}^{t} (Ax(s)+f(s))ds$ for all $t\in I$  and $y(t)=Hx(t)+\eta(t)$ a.e. on $I$. Denote the set of all solutions on $I$ by $\mathcal{B}_I(F,A,H)$.
\end{Definition}
Let us fix the symmetric positive definite matrices $Q_0,Q,R$.
\begin{Definition}
A \emph{solution} $(x,f,y,\eta) \in \mathcal{B}_{I}(F,A,H)$ is said to be admissible if $\rho(Fx(0),f,\eta,t_1,Q_0,Q,R) \le 1$ for all $t_1 \in I$.
Denote the set of all admissible solutions by $\mathscr{EE}(I)$.
\end{Definition}
\begin{Remark}\label{r:rho_norm_Fx}
The map $(x_0,f,\eta)\mapsto\rho(x_0,f,\eta,t_1,Q_0,Q,R)$ is a norm of the space $\mathcal H_{t_1}:=\mathbb{R}^m \times L^2(I,\mathbb{R}^m)  \times  L^2(I,\mathbb{R}^p)$ and
$\mathscr{ E}(t_1):=\{(x_0,f,\eta):\rho(x_0,f,\eta,t_1,Q_0,Q,R)\le 1\}$
is the unit ball of $\mathcal H_{t_1}$ induced by this norm. 
\end{Remark}

Next, we define finite horizon minimax observers.
Fix a time $0 < t_1 \in \mathbb{R}$ and set $I=[0,t_1]$. Fix $\ell\in\mathbb R^m$.
To define an observer we rely upon classical results~\cite{KrenerMinimax,Tempo1985,Nakonechnii1978} on minimax (worst-case) estimators (observers) which state that for general linear equations the minimax observer is linear in outputs $y$~\cite{Tempo1985}. The latter suggests that an estimate $\widehat{\ell^{\top}Fx(t_1)}$ of the true value $\ell^{\top}Fx(t_1)$, which is computed by the minimax observer, should be linear in $y$. It is then natural to model the effect of applying minimax observers to outputs as continuous linear functionals, generated by functions $U \in L^2(I,\mathbb{R}^p)$ as follows:
\[ \obs{t_1}:=\int_{0}^{t_1} y^{\top}(t)U(t)dt\,, \quad \forall y \in L^2(I,\mathbb{R}^p)\,. \]
Hence, given $y$, we compute the estimate $\widehat{\ell^{\top}Fx(t_1)}$ of the true value $\ell^{\top}Fx(t_1)$ by selecting a function $U$
and evaluating $\obs{t_1}$. In order to select $U$ optimally we introduce a cost function $\sigma$ measuring the worst-case estimation error:
\begin{equation}
  \label{eq:error:fin}
    \sigma(U,\ell,t_1):=\sup_{(x,f,y,\eta)\in \mathscr{EE}(I)}(\ell^{\top}Fx(t_1)-\obs{t_1})^2\,.
\end{equation}
\begin{Definition}\label{d:finhorobs}
We say that $\widehat{U}\in L^2(I,\mathbb{R}^p)$ is a \emph{minimax observer on $I=[0,t_1]$}, if \( \hat\sigma:=\sigma(\widehat{U},\ell,t_1)=\inf_{U \in L^2(I,\mathbb{R}^p)} \sigma(U,\ell,t_1) < +\infty \). 
\end{Definition}
If $\sigma(U,\ell,t_1) < +\infty$, then the estimation error satisfies $(\ell^{\top}Fx(t_1) - \mathscr{O}_{U}(y)(_1t))^2  \le \sigma(U,\ell,t_1)$ for \textbf{any} admissible solution $(x,f,y,\eta) \in \mathscr{EE}(I)$.  If, in addition, $U$ is a minimax observer on $I$, then the error bound $\sigma(U,\ell,t_1)$ is the
smallest possible.
We show below that if there exists a function $U \in L^2(I,\mathbb{R}^p)$ with finite worst-case error, then there also exists the minimax  observer $\widehat U$,
and it can be implemented as an output of a linear system.

Consider now the infinite horizon case, i.e. $I=[0,+\infty)$.
Similarly to the finite horizon case, observers of interest are devices which map outputs $y$ to estimates $\widehat{\ell^{\top}Fx(t)}$ of $\ell^{\top}Fx(t)$, and 
 $\widehat{\ell^{\top}Fx(t)}$  is a linear and continuous function of $y|_{[0,t)}$.
To model such observers we define $\mathscr{F}$ -- the set of all maps
  $U:\{ (\tau,s) \mid \tau > 0, s \in [0,\tau]\} \rightarrow \mathbb{R}^p$ such that $\forall\tau > 0$,
  the function $U(\tau,\cdot): [0,\tau] \ni s \mapsto U(\tau,s)$ belongs to $L^2([0,\tau],\mathbb{R}^p)$. An element $U \in \mathscr{F}$ maps output $y \in L^2_{loc}(I,\mathbb{R}^p)$ to a function $\mathscr{O}_{U}(y)$ defined as:
 \[ \mathscr{O}_{U}(y): I\ni t\mapsto \mathscr{O}_U(y)(t):=\int_0^{t} y^{\top}(s)U(t,s)ds\,.\]
Now, to construct $\widehat{\ell^{\top}Fx(t)}$ we take $U \in \mathscr{F}$ and compute $\mathscr{O}_U(y)$. To select $U$ optimally we define the worst-case estimation error:
\begin{equation}
\label{problem:obs:eq2.1}
  \sigma(U,\ell):=\limsup_{t \rightarrow \infty}\sigma(U(t,\cdot),\ell,t)\,.
\end{equation}
\begin{Definition}
\label{d:infhorobs}
We say that $\widecheck{U} \in \mathscr{F}$ is an \emph{infinite horizon minimax observer}, if
   \begin{equation}
   \label{problem:obs:eq2}
     \sigma(\widecheck{U},\ell)= \inf_{U\in\mathscr{F}} \sigma(U,\ell) < +\infty.
   \end{equation}
\end{Definition}
If $\sigma(U,\ell) < +\infty$, then for any $\epsilon > 0$, for large enough $t$, the estimation error satisfies $(\ell^{\top}Fx(t) - \mathscr{O}_{U}(y)(t))^2  \le \sigma(U,\ell)+\epsilon$ for \textbf{any} admissible solution $(x,f,y,\eta) \in \mathscr{EE}([0,+\infty))$.  If, in addition, $U$ is an infinite horizon minimax observer, then
the error bound $\sigma(U,\ell)$ is the smallest possible.
We will be interested in infinite horizon minimax observers which can be represented by stable LTI systems.
\begin{Definition}
$U \in \mathscr{F}$ is said to be represented by an LTI system $(A_o,B_o,C_o,0)$,  
if for any $y \in L^2_{loc}(I,\mathbb{R}^p)$, where
$I=[0,t_1]$, $0 < t_1 \in \mathbb{R}$ or $I=[0,+\infty)$,
 it holds that
 $\forall t\in I:\mathscr{O}_{U}(y)(t)=C_os(t)$, where $\dot s(t) = A_os(t)+B_oy(t)$, $s(0)=0$. 
\end{Definition}
 Note that $(A_o,B_o,C_o,0)$ represents an element $U$ of $\mathscr{F}$, if and only if $U(t,s)=B_o^{\top}e^{A^{\top}_o(t-s)}C_o^{\top}$ for all $t > 0, s \in [0,t]$.
 Moreover, we will prove that if there exist a $U \in \mathscr{F}$ with finite worst-case estimation error, then there exist an infinite horizon minimax observer $\widecheck{U}$, which, in addition, can be represented by a stable LTI system.


%
\section{Main results}
\label{sec:main}
\subsection{Duality and existence of minimax observers}
\label{sec:dual}
We show that the minimax observer is related to the solution of an optimal control problem for the dual DAE:
\begin{align}
&\frac{d(F^{\top}q(t))}{dt}=A^{\top}q(t)-H^{\top}u(t)\,.   
\label{dae:dual}
\end{align}
To this end, we need to define what we mean by a solution of \eqref{dae:dual}.
\begin{Definition}
Let $I$ be an interval of the form $I=[0,t_1]$, $0<t_1 \in \mathbb{R}$ or $I=[0,+\infty)$.
A pair $(q,u) \in L^2_{loc}(I,\mathbb{R}^m) \times L^2_{loc}(I,\mathbb{R}^p)$
is a solution of \eqref{dae:dual} on $I$, if $F^{\top}q$ is absolutely continuous and $F^{\top}q(t)=F^{\top}q(0)+\int_0^t (A^{\top}q(s)-H^{\top}u(s))ds$ for all $t \in I$. Denote the set of all solutions on $I$ by $\mathcal{B}_I(F^{\top},A^{\top},H^{\top})$ and define \( \mathscr{D}_{\ell}(I):=\{ (q,u) \in \mathcal{B}_I(F^{\top},A^{\top},H^{\top})\mid F^{\top}q(0)=F^{\top}\ell\}\).
\end{Definition}
\begin{Definition}[$\ell$-impulse observability]
\label{d:l-obs}
The tuple $(F,A,H)$ is \emph{$\ell$-impulse observable}, if $\mathscr{D}_{\ell}([0,t_1])\ne\varnothing$ for some $0<t_1 \in \mathbb{R}$. 
\end{Definition}
We define now the following cost function for \eqref{dae:dual}. Fix $0 <  t_ 1 \in \mathbb{R}$ and let $I$ be an interval of the form $I=[0,\tau]$ for some $t_1 \le \tau \in \mathbb{R}$, or $I=[0,+\infty)$.
For any $(q,u) \in \mathscr{D}_{\ell}(I)$, define
\begin{align}
   J(q,u,t_1):=\rho(q(t_1),q,u,t_1,\bar Q_0,Q^{-1},R^{-1})\,,\label{control:dual}
\end{align}
where \(\bar{Q}_0:=F({F^{\top}}^+ -\LOPT)^{\top}Q_0^{-1} ({F^{\top}}^+-\LOPT)F^{\top}\) with \(\LOPT=P(PQ_0^{-1}P)^+PQ_0^{-1}{F^{\top}}^{+}\) and \(P=(\Id_m-(F^{\top})^+F^{\top})\) and $\rho$ is defined in ~\eqref{eq:rhox0feta} with
$\tau=t_1$, $g_0=q(t_1)$,  $g=q$, $h=u$, $S_0=\bar{Q}_0$, $S_1=Q^{-1}$, $S_2=R^{-1}$.
We can now state the following theorem relating existence of minimax observers to minimizing the cost function \eqref{control:dual}.
\begin{Theorem}\label{p:1:theo}
The following statements are equivalent:
\begin{itemize}
\item [(i)] $(F,A,H)$ is $\ell$-impulse observable,
\item [(ii)] there exists $U\in L^2([0,t_1],\mathbb{R}^p)$ such that $\sigma(U,\ell,t_1)<+\infty$,
 \item [(iii)] there exists a minimax observer $\widehat U$ on $[0,t_1]$, such that  $\widehat{U}=\delta_{t_1}(u^*)$ and \( \sigma(\widehat{U},\ell,t_1)=J^*_{t_1}\), where $(q^*,u^*) \in \mathscr{D}_{\ell}([0,t_1])$ is a solution of the dual control problem: \( J^*_{t_1}: = J(q^*,u^*,t_1)=\inf_{(q,u) \in \mathscr{D}_{\ell}(I)} J(q,u,t_1) \).
\end{itemize}
Moreover, for any $U\in L^2([0,t_1],\mathbb{R}^p)$,
\begin{equation}
   \label{eq:sigmaUt1}
    \sigma(U,\ell,t_1)=\inf_{(q,u) \in \mathscr{D}_{\ell}(I), u=\delta_{t_1}(U)} J(q,u,t_1)\,.
 \end{equation}
In particular, the set $\{(q,u) \in \mathscr{D}_{\ell}(I) \mid u=\delta_{t_1}(U)\}$ is not empty, if and only if $\sigma(U,\ell,t_1)<+\infty$.
\end{Theorem}
Next, we present the  notion of $\ell$-detectability and a new duality principle for the infinite horizon case.
\begin{Definition}[$\ell$-detectability]
\label{d:l-detect}
The tuple $(F,A,H)$ is \emph{$\ell$-detectable}, if there exists $(q,u) \in \mathscr{D}_{\ell}([0,+\infty))$ such that $\lim_{t \rightarrow \infty} F^{\top}q(t)=0$.
\end{Definition}
  Note that definition~\ref{d:l-detect} implies that $\mathscr{D}_{\ell}([0,t_1])\ne\varnothing$ for some $0<t_1 \in \mathbb{R}$: indeed, for $(q,u)$ from the definition~\ref{d:l-detect} it follows that its restriction $(q|_{[0,t_1]},u|_{[0,t_1]})$ onto $[0,t_1]$ belongs to
$ \mathscr{D}_{\ell}([0,t_1])$ and hence for any $t_1 > 0$
$ \mathscr{D}_{\ell}([0,t_1]) \ne \emptyset$. Hence, $\ell$-detectability implies $\ell$-impulse observability by definition~\ref{d:l-obs}.
\begin{Theorem}\label{p:5} The following statements are equivalent:
\begin{itemize}
\item[(i)] $(F,A,H)$ is $\ell$-detectable,
\item[(ii)] there exists an infinite horizon minimax observer $\widecheck{U}$, such that $\widecheck{U}(t_1,s)=\delta_{t_1}(u^{*})(s)$, $0 \le s \le t$ and $\sigma(\widecheck{U},\ell)=J^*_\infty$, where $(q^{*},u^{*}) \in \mathscr{D}_{\ell}([0,+\infty)$ satisfies 
\begin{equation}
  \label{eq:dual_LQ_inf}
 \begin{split}
 & J^*_\infty:=\limsup_{\tau \rightarrow \infty } J(q^{*},u^{*},\tau)\\
  & = \limsup_{\tau \rightarrow \infty } \inf_{(q,u) \in \mathscr{D}_{\ell}([0,\tau])} J(q,u,\tau) < +\infty\,.
  \end{split}
\end{equation}
\end{itemize}
\end{Theorem}
 The assumption $Q_0,Q,R >0$ is essential for Theorems ~\ref{p:1:theo} -- ~\ref{p:5}.
\begin{Remark}
\label{rem_impulse_obs}
Note that $(F,A,H)$ is $\ell$-impulse observable \emph{for all $\ell \in \mathbb{R}^m$}, if and only if
it is impulse controllable in the sense of  \cite{DAEBookChapter}. The latter can be characterized in terms of $(F,A,H)$.
If $\mathrm{rank} \begin{bmatrix} F^{\top} & A^{\top} & H^{\top} \end{bmatrix}^{\top} = n$, then $(F,A,H)$ is $\ell$-impulse observable for all $\ell \in \mathbb{R}^m$ if and only if $(F,A,H)$ is impulse observable in the sense of \cite{Darouach2009,HouMuelle1999,IshiharaTerra}. 
The latter is equivalent to
$\rank~ \left[\begin{smallmatrix} F & A \\ 0 & H \\ 0 &  F \end{smallmatrix}\right] = n + \rank~ F$. 
We stress that the condition above characterizes when it is true that for all $\ell$, $(F,A,H)$ is
$\ell$-impulse observable. Hence, this condition does not involve $\ell$. In contrast,
the conditions which characterize $\ell$-impulse observability of $(F,A,H)$ for a \emph{given} $\ell$ will depend on that vector $\ell$. 
%
%
 From ~\cite[Lemma 4.3.7]{BergerThesis} it follows that $(F,A,H)$ is $\ell$-detectable for all those $\ell \in \mathbb{R}^m$ for which it is $\ell$-impulse observable $\Leftrightarrow$ 
 $\forall \lambda \in \mathbb{C}, \mathrm{Re} \lambda \ge 0: \mathrm{rank} \left[\begin{smallmatrix} \lambda F - A \\  H \end{smallmatrix}\right] = \max_{s \in \mathbb{C}} \mathrm{rank}\left[ \begin{smallmatrix} sF - A \\ H \end{smallmatrix}\right]$.
\end{Remark}

\subsection{Observer design}
\label{sec:solut-dual-contr}
Below we present algorithms for computing minimax observers. Recall from \cite{AutomaticaPaper} the notion of LTI systems associated with \eqref{dae:dual}.
\begin{Definition}
\label{linassc:def}
  Let $A_a \in \mathbb{R}^{\hat{n} \times \hat{n}}$, $B_a \in \mathbb{R}^{\hat{n} \times k}$, $C_a \in \mathbb{R}^{(m+p) \times \hat{n}}$, $D_a \in \mathbb{R}^{(m+p) \times k}$ for some $\hat{n},k > 0$.  The
 LTI system $\mathscr{S}=(A_a,B_a,C_a,D_a)$ is said to be
\emph{associated with ~\eqref{dae:dual}}, if the following conditions hold:
  \begin{enumerate}
\item Either $k=1$ and $\left[\begin{smallmatrix} D_a\\ B_a \end{smallmatrix}\right]=0$ or $D_a$ is full column rank. 
\item
  Let $C_s$ and $D_s$ be the matrices formed by the first
  $m$ rows of $C_a$ and $D_a$ respectively. Then $F^{\top}D_s=0$ and $\Rank~ F^{\top}C_s = \hat{n}$.
 \item
   Let $I$ be an interval of the form $I=[0,t_1]$, $0<t_1 \in \mathbb{R}$ or $I=[0,+\infty)$.
    Then $(q,u) \in \mathcal{B}_I(F^{\top},A^{\top},H^{\top})$ if and only if there exists
       $(p,g)  \in AC(I,\mathbb{R}^{\hat{n}}) \times L^2_{loc}(I,\mathbb{R}^k)$
       such that $\dot p =A_ap+B_ag$ a.e. and $
       \left[\begin{smallmatrix}
         q\\u
       \end{smallmatrix}\right]=C_ap+D_ag$ a.e.
 \end{enumerate}
The matrix $\LMAP:=(F^TC_s)^+$ is called the state map of $\mathscr{S}$.  
\end{Definition}
Associated LTI systems allow us to represent solutions of DAEs as outputs of linear systems.
Recall from \cite[Theorem 2-3]{AutomaticaPaper} that there exists an LTI system $\mathscr{S}=(A_a,B_a,C_a,D_a)$ associated with \eqref{dae:dual} and that all LTI systems which are associated with \eqref{dae:dual} are feedback equivalent.\\ The matrix $\LMAP$ allows us to relate the state trajectories of \eqref{dae:dual} and the state trajectories of $\mathscr{S}$:
 from \cite[Theorem 1]{AutomaticaPaper} it follows that  for any $(q,u)\in \mathcal{B}_I(F^{\top},A^{\top},H^{\top})$,  for $v:= \LMAP F^{\top}q$ and $g:=D_a^{+}(\left[\begin{smallmatrix}
         q\\u
       \end{smallmatrix}\right] -C_av)$, $\dot v = A_av+B_ag$ a.e. and $\left[\begin{smallmatrix} q\\u       \end{smallmatrix}\right]  = C_av + D_ag$ a.e.
Consequently (see \cite[Corollary 1]{AutomaticaPaper}) it follows that $(F,A,H)$ is $\ell$-impulse observable if and only if $F^\top \ell \in \mathrm{Im} F^\top C_s$. An algorithm for computing the associated LTI was discussed in \cite[Proof of Th.1]{AutomaticaPaper} and it is summarized in Algorithm \ref{alg1}. This algorithm was implemented, the Matlab function $DAE2Lin$ implementing Algorithm \ref{alg1} is available in the file
$main\_TA3.m$ of the supplemenatry material of this report.
\begin{algorithm}
\caption{Computation of the associated LTI system ~\cite[Proof of Th.1]{AutomaticaPaper}. \label{alg1}}
\renewcommand{\algorithmicrequire}{\textbf{Input:}}
\renewcommand{\algorithmicensure}{\textbf{Output:}}
\begin{algorithmic}[1]
\REQUIRE $(F,A,H)$
\ENSURE LTI system $\mathscr{S}=(A_a,B_a,C_a,D_a)$ associated with the dual system \eqref{dae:dual} and the state map $\LMAP$.

\STATE Compute the SVD decomposition
$F^T=U\Sigma V^T$ of $F^T$ and assume that
$\Sigma=\diag(\sigma_1,\ldots,\sigma_r, 0,\ldots, 0) \in \mathbb{R}^{n \times m}$, $r=\rank ~ F$,
$U \in \mathbb{R}^{n \times n}, V \in \mathbb{R}^{m \times m}$, $U^{T}U=\mathcal{I}_n$, $V^TV=\mathcal{I}_m$.
       \[
          \begin{split}
          S=\diag(\frac{1}{\sqrt{\sigma_1}},\ldots, \frac{1}{\sqrt{\sigma_r}},\underbrace{1,\ldots,1}_{n-r}) U^T, \\
          T=V\diag(\frac{1}{\sqrt{\sigma_1}},\ldots, \frac{1}{\sqrt{\sigma_r}},\underbrace{1,\ldots,1}_{m-r}).
         \end{split}
       \]
Note that
    \( SF^TT = \begin{bmatrix} I_r & 0 \\
                          0   & 0
          \end{bmatrix} \).

\STATE
 Let $\widetilde{A} \in \mathbb{R}^{r \times r}$, $B_{1} \in \mathbb{R}^{r \times p}$, $G \in \mathbb{R}^{r \times \kappa}$, $\widetilde{D} \in \mathbb{R}^{(n-r) \times \kappa}$,
$\kappa=m-r+p$, 
$\widetilde{C} \in \mathbb{R}^{(n-r) \times r}$,  be such that
 \[
    \begin{split}
    & SA^TT=\begin{bmatrix} \widetilde{A} & A_{12} \\
                    A_{21} & A_{22}
     \end{bmatrix}, ~  SH^T=-\begin{bmatrix} B_1 \\ B_2 \end{bmatrix} \\
    & G = \begin{bmatrix} A_{12}, & B_1 \end{bmatrix}  \mbox{, \ \ }
     \widetilde{D} = \begin{bmatrix} A_{22}, & B_2 \end{bmatrix}  \mbox{ and }
     \widetilde{C}=A_{21}.
   \end{split}
 \]

\STATE
   Consider the linear system
   \begin{equation}
  \label{eq:DAE_lin}
    \begin{split}
     & \dot p = \widetilde{A}p+Gq, ~  z     = \widetilde{C}p + \widetilde{D}q.
    \end{split}
 \end{equation}

  Using  the algorithms discussed in \cite[Section 4.5, page 241]{GeomBook2}, compute a full row rank matrix $\mathscr{V}$ and a matrix $\widetilde{F} \in \mathbb{R}^{\kappa \times r}$, where $\kappa$ is the number of columns of $G$,
  such that
  \begin{itemize}
  \item
   $\IM \mathscr{V}$ is   the largest output nulling controlled invariant subspace of \eqref{eq:DAE_lin}, and
  \item
  $(\widetilde{A}+G\widetilde{F})(\IM \mathscr{V}) \subseteq \IM \mathscr{V}$.
  \end{itemize}
  For instance, $(\mathscr{V},\widetilde{F})$ can be computed by applying the Matlab function $vstar$  from \cite{MarroAlg} to $(\widetilde{A},G,\widetilde{C},\widetilde{D})$.

\STATE
Let $L \in \mathbb{R}^{\kappa \times k}$ be such that $\IM L=\ker \widetilde{D} \cap G^{-1}(\V)$, and either $k=\Rank L$, or $L=(0,\ldots,0)^T \in \mathbb{R}^\kappa$.

For instance, using Matlab functions $invt$,$ints$,$ker$ of \cite{MarroAlg}, the function call   $ints(invt(G,\mathscr{V}),ker(\widetilde{D}))$ will return $L$.

\STATE
   Define
   \[
   \begin{split}
    &   \bar{C} = \begin{bmatrix} T & 0 \\ 0 & I_n \end{bmatrix} \begin{bmatrix} I_r \\ \widetilde{F} \end{bmatrix} \mbox{ and }
       \bar{D} = \begin{bmatrix} T & 0 \\ 0 & I_n \end{bmatrix}\begin{bmatrix} 0 \\ L \end{bmatrix}, \\
   & A_a=P_{\mathscr{V}} (\widetilde{A}+G\widetilde{F})\mathscr{V}, B_a= P_{\mathscr{V}} GL,  \quad C_a =  \bar{C}  \mathscr{V}, \\
   & D_a=\bar{D}, ~  \LMAP=(F^{\top}C_s)^{+}, 
  \end{split}
 \]
 where $P_{\mathscr{V}}=\mathscr{V}^{+}$, and
 $C_s$ is the matrix formed by the first $m$ rows of $C_a$.


\end{algorithmic}
\end{algorithm}
Now we describe how to construct a minimax observer $\widehat{U}$ on the interval $I=[0,t_1]$, $t_1  >0$.
Let $\mathscr{S}=(A_a,B_a,C_a,D_a)$ be an LTI system associated with \eqref{dae:dual}. Consider the equation:
  \begin{equation}
  \label{DARE}
  \begin{split}
     & \dot P(t)\!\! = \!\! A_a^{\top}P(t)+P(t)A_a\!\!-\!\!K^{\top}(t)(D_a^{\top}SD_a)K(t)+C_a^{\top}SC_a\,, \\
     & P(0)=(F^{\top}C_s)^{\top}\bar{Q}_0F^{\top}C_s, \quad
 S=\left[\begin{smallmatrix} Q^{-1},& 0 \\ 0, & R^{-1} \end{smallmatrix}\right]\,, \\
     & K(t)\!\!=\!\! \left\{\!\!\! \begin{array}{rl}
                (D_a^{\top}SD_a)^{-1}(B_a^{\top}P(t)+D_a^{\top}SC_a) & \!\!\! \mbox{ if }  D_a\ne 0\,, \\
                0  \in \mathbb{R}^{1 \times \hat{n}}  & \!\!\! \mbox{ if }  D_a = 0\,.
             \end{array}\right.
  \end{split}
  \end{equation}
   where $C_s$ is as in Definition \ref{linassc:def}. Note that \eqref{DARE} is well defined, as either $D_a$ is full column rank, and so  $D_a^{\top}SD_a$ is invertible, or $B_a=0,D_a=0,k=1$.
 Let $u^* \in L^2(I,\mathbb{R}^p)$, $v^* \in AC(I,\mathbb{R}^{\hat{n}})$ be such that
 \begin{equation}
 \label{dual:finhor:sol}
 \begin{split}
   & u^*(s)=(C_{u}-D_{u}K(t_1-s))v(s)\,, \\
   & \dot v(s) = (A_a-B_aK(t_1-s))v(s) \mbox{,  \ \ } v(0)=\LMAP(F^{\top}\ell)\,,
 \end{split}
 \end{equation}
 where $C_u$ and $D_u$ are the matrices formed by the last $p$ rows of $C_a$ and $D_a$ respectively, and
 $K$ is as in \eqref{DARE}.
\begin{Theorem}[Minimax observer: finite horizon]
\label{t:observer-fin-hor}
If $(F,A,H)$ is $\ell$-impulse observable, then $\widehat{U}(s)=u^{*}(t_1-s)$ is a finite horizon minimax observer on $I=[0,t_1]$. Moreover, \( \sigma(\widehat{U},\ell,t_1)=\ell^{\top}F\LMAP^{\top} P(t_1)\LMAP F^{\top}\ell \)
  and $\mathcal{O}_{\widehat{U}}(t_1,y) = \ell^{\top}F\LMAP^{\top}r(t_1)$, where $r \in AC(I,\mathbb{R}^{\hat{n}})$, $r(0)=0$ and
  \begin{equation}
  \label{finhor:eq1}
     \dot r(t)=(A_a-B_a K(t))^{\top}r(t)+(C_u-D_uK(t))^{\top}y(t) \quad \mbox{ a.e.}
   \end{equation}
\end{Theorem}
Note that $\ell$-impulse observability can easily be checked, see Algorithm \ref{alg5} below. The Matlab function $IsObservable$ implementing Algorithm \ref{alg5} is available
in the file $main\_TA3.m$ of supplementary material of this report.
\begin{algorithm}[H]
\caption{Checking $\ell$-impulse observability \label{alg5}}
\renewcommand{\algorithmicrequire}{\textbf{Input:}}
\renewcommand{\algorithmicensure}{\textbf{Output:}}
\begin{algorithmic}[1]
\REQUIRE $(F,A,H,\ell)$
\ENSURE Yes, if $(F,A,H)$ is $\ell$-impulse observable, No otherwise
\STATE Use Algorithm \ref{alg1} to compute an LTI system $\mathscr{S}=(A_a,B_a,C_a,D_a)$ associated with \eqref{dae:dual}.
\STATE
\label{step:check2}
Let $C_s$ be the matrix formed by the first $m$ rows of $C_a$. Check if $\ell$ belongs to $\IM (F^TC_s)$, for example, by checking if $F^TC_s(F^TC_s)^{+}\ell=\ell$.
\STATE If $\ell \in \IM (F^TC_s)$, return Yes, otherwise return No.
\end{algorithmic}
\end{algorithm}
Theorem \ref{t:observer-fin-hor} yields an algorithm for computing the finite horizon minimax estimate. This algorithm is presented in  Algorithm \ref{alg2} below.
\begin{algorithm}[H]
\caption{Finite horizon minimax state estimate $\mathcal{O}_{\widehat{U}}(t_1,y)$ \label{alg2}}
\renewcommand{\algorithmicrequire}{\textbf{Input:}}
\renewcommand{\algorithmicensure}{\textbf{Output:}}
\begin{algorithmic}[1]
\REQUIRE $(F,A,H,\ell)$ and the observed output signal $y \in L^2([0,t_1],\mathbb{R}^p)$
\ENSURE Minimax estimate $\mathcal{O}_{\widehat{U}}(t_1,y)$
\STATE Use Algorithm \ref{alg1} to compute the linear system $\mathscr{S}=(A_a,B_a,C_a,D_a)$ associated with \eqref{dae:dual} and the state map $\LMAP$.
\STATE Compute $P$ and $K$ by solving~\eqref{DARE} and find $r(t)$ by solving~\eqref{finhor:eq1}. The numerical approximations of $P$ and $r$ may be done by using M\"{o}bius time integrators~\cite{ZhukCDC14}.
\STATE Return $\mathcal{O}_{\widehat{U}}(t_1,y) = \ell^{\top}F\LMAP^{\top}r(t_1)$.
\end{algorithmic}
\end{algorithm}
 Next, we present an algorithm for computing an infinite horizon minimax observer $\widehat{U}$.
 To this end, recall from \cite{AutomaticaPaper} the notion of a stabilizable LTI system associated with the DAE~\eqref{dae:dual}.
 Let $\mathscr{S}=(A_a,B_a,C_a,D_a)$ be an LTI system associated with \eqref{dae:dual} and consider the stabilizability subspace $\V_g$ of $\mathscr{S}$. From~\cite{CalierDesoer} it then follows that $\V_g$ is $A_a$-invariant and $\IM B_a \subseteq \V_g$. Hence, there exists a basis transformation $T$ such that $T(\V_g)=\IM \left[\begin{smallmatrix} I_l & 0 \\ 0 & 0 \end{smallmatrix}\right]$, $l=\dim \V_g$ and in this new basis,
  \(TA_aT^{-1}=\left[\begin{smallmatrix} A_g & \star \\ 0 & \star \end{smallmatrix}\right]\), $
      TB_a=\left[\begin{smallmatrix} B_g \\ 0 \end{smallmatrix}\right]$,
       $C_aT^{-1}=\left[\begin{smallmatrix} C_g^{\top} \\ \star \end{smallmatrix}\right]^{\top}$, where
  $A_g \in \mathbb{R}^{l \times l}$, $B_g \in \mathbb{R}^{l \times k}$, $C_g \in \mathbb{R}^{(n+m) \times l}$.
\begin{Definition}
\label{assco_stab:def}
  The LTI system $\mathscr{S}_g=(A_g,B_g,C_g,D_g)$ with $D_g=D_a$ is called a stabilizable LTI system associated with \eqref{dae:dual}.  The matrix $\LMAP_g=\left[\begin{smallmatrix} I_l & 0 \end{smallmatrix}\right] T\LMAP$ is called the associated state map.
 \end{Definition}
The LTI system $\mathscr{S}_g$ is a restriction of $\mathscr{S}$ to the subspace $\V_g$. It follows that $\mathscr{S}_g$ is stabilizable.
Moreover, since all associated LTI systems of \eqref{dae:dual} are feedback equivalent, all stabilizable associated LTI systems of~\eqref{dae:dual} are feedback equivalent.
From the discussion after Lemma 1 of \cite{AutomaticaPaper} it follows that the LTI $\mathscr{S}_g$ and the map $\LMAP_g$ can be computed from the matrices $F,A,H$.
For the convenience of the reader, we present the algorithm in  Algorithm~\ref{alg3} below. The Matlab function $Lin2StabLin$ implementing Algorithm \ref{alg3} can be found in 
the file $main\_TA3.m$ in supplementary material of this report.
\begin{algorithm}[H]
\caption{Stabilizable associated linear system
~\cite[Lem.1 and discussion after]{AutomaticaPaper}.
\label{alg3}}
\renewcommand{\algorithmicrequire}{\textbf{Input:}}
\renewcommand{\algorithmicensure}{\textbf{Output:}}
\begin{algorithmic}[1]

\REQUIRE $(F,A,H)$
\ENSURE Stabilizable LTI system $\mathscr{S}_g=(A_g,B_g,C_g,D_g)$ associated with \eqref{dae:dual} and the state map $\LMAP_g$.

\STATE Use Algorithm \ref{alg1} to compute the LTI system $\mathscr{S}=(A_a,B_a,C_a,D_a)$ associated with \eqref{dae:dual} and the state map $\LMAP$.

\STATE
\label{step:stable1}
       Compute the controllability matrix $R$ of $\mathscr{S}$,
\STATE
       Compute a matrix $\mathscr{V}_2$  such that
       $\IM \mathscr{V}_2$ is largest $A_a$-invariant subspace such that the restriction of $A_a$ to $\IM \mathscr{V}_2$ is stable.
       If we use the Matlab tool \cite{MarroAlg}, then $\mathscr{V}_2$ can be taken as the first component of the output of the $subsplit(A_a)$.
\STATE
\label{step:stable2}
       Let $\mathscr{V}_g$ be any matrix whose columns form a basis of $\IM \begin{bmatrix} R, & \mathscr{V}_2 \end{bmatrix}$. Then $\mathscr{V}$ is full column rank, and
       from \cite[Theorem 4.30, page 96]{TrentelmanBook} it follows that $\IM \mathscr{V}_g$ equals the largest stabilizable subspace $\V_g$ of $\mathscr{S}$.

   In Matlab, $\mathscr{V}_g$ can be taken as the output of $orth(\begin{bmatrix} R,\mathscr{V}_2 \end{bmatrix})$.

\STATE
       \( A_g=P_{\mathscr{V}_g} A_a \mathscr{V}_g, ~ B_g=P_{\mathscr{V}_g} B_a, ~ C_a=C_g\mathscr{V}_g, ~ \LMAP_g = P_{\mathscr{V}_g} \LMAP,  \)
      where $P_{\mathscr{V}_g}=\mathscr{V}_g^{+}$.
\end{algorithmic}
\end{algorithm}
Note that by \cite{AutomaticaPaper}, $(F,A,H)$ is $\ell$-detectable if and only if $F^T\ell \in \IM (F^TC_s)$ and $\LMAP(F^T\ell)$ belongs to the stabilizability subspace $\V_g$ of $\mathscr{S}$.
This remark yields the following algorithm, presented in Algorithm \ref{alg6}, for checking $\ell$-detectability.
\begin{algorithm}[H]
\caption{Checking $\ell$-detectability \label{alg6}}
\renewcommand{\algorithmicrequire}{\textbf{Input:}}
\renewcommand{\algorithmicensure}{\textbf{Output:}}
\begin{algorithmic}[1]
\REQUIRE $(F,A,H,\ell)$
\ENSURE Yes, if $(F,A,H)$ is $\ell$-detectable, No otherwise
\STATE Use Algorithm \ref{alg2} to compute an LTI system $\mathscr{S}=(A_a,B_a,C_a,D_a)$ associated with \eqref{dae:dual} and the state-map $\LMAP$.
\STATE Use Steps \ref{step:stable1}-\ref{step:stable2} of Algorithm \ref{alg3} to compute a full column rank matrix $\mathscr{V}_g$ such that $\IM \mathscr{V}_g$ is the largest stabilizable subspace $\V_g$ of $\mathscr{S}$.

\STATE Like in Step \ref{step:check2} of Algorithm~\ref{alg5}, check if $\ell \in \IM F^TC_s$, where $C_s$ is the matrix formed by the first $m$ rows of $C_a$. If $\ell \notin \IM F^TC_s$, then quit the algorithm and return No.
\STATE If $\ell \in \IM F^TC_s$, then check if $\LMAP(F^T\ell) \in \IM \mathscr{V}_g$. The latter can be done, for example, by checking if
 $\mathscr{V}_g\mathscr{V}_g^{+}(\LMAP(F^T\ell))=\LMAP(F^T\ell)$.

\STATE If $\LMAP(F^T\ell) \in  \IM \mathscr{V}_g$, then return Yes, otherwise return No.
\end{algorithmic}
\end{algorithm}
The Matlab function $IsDetectable$  implementing Algorithm \ref{alg6} can be found in the file $main\_TA3.m$
of the supplementary material of this report.

Let $\mathscr{S}_g=(A_g,B_g,C_g,D_g)$ and $\LMAP_g$ be as in Definition \ref{assco_stab:def}. 
Consider the following algebraic Riccati equation:
    \begin{equation}
    \label{ARE}
    \begin{split}
       & 0 = PA_g+A^{\top}_gP- K^{\top}(D^{\top}_gSD_g)K + C^{\top}_gSC_g\,, \\
       & K=\left\{ \begin{array}{rl}
          (D^{\top}_gSD_g)^{-1}(B^{\top}_gP+D^{\top}_gSC_g)      & \mbox{ if } D_g \ne 0\,, \\
         0 \in \mathbb{R}^{1 \times l} & \mbox{ if } D_g=0\,,
       \end{array}\right. \\
     & S=\left[\begin{smallmatrix} Q^{-1} & 0 \\ 0 & R^{-1} \end{smallmatrix}\right]\,.
    \end{split}
    \end{equation}
Note that either $D_g=D_a$ is full column rank, and so, by positive definiteness of $Q,R$, $D_g^{\top}SD_g$ is invertible, or $B_g=0,D_g=0,k=1$.
%
If $P=P^T > 0$ is a solution of \eqref{ARE}, then  let $u^* \in L^2_{loc}([0,+\infty),\mathbb{R}^p)$, $v^* \in AC([0,+\infty),\mathbb{R}^l)$ be such that
 \begin{equation}
      \label{opt8}
        \begin{split}
        & u^*(t)=(\hat{C}_u-\hat{D}_uK)v^*(t)\,,\\
        & \dot {v^*}(t) = (A_g-B_gK)v^*(t)\,, \quad v^*(0)=\LMAP_g(F^{\top}\ell)\,,
        \end{split}
 \end{equation}
 where $\hat{C}_u$ and $\hat{D}_u$ are the matrices formed by the last $p$ rows of $C_g$ and $D_g$ respectively. Define
 \[
     A_o=(A_g-B_gK)^{\top}, \quad B_o=(\hat{C}_u-\hat{D}_uK)^{\top}, \quad C_o=\ell^{\top}F\LMAP^T_g\,. \]
 \begin{Theorem}[Minimax observer: infinite horizon]\label{t:observer-infin-hor}
 \label{main:obs}
  Assume that $(F,A,H)$ is $\ell$-detectable. Then the following holds:
  \begin{itemize}
  \item [(i)] $A_o$ is a Hurwitz matrix and \eqref{ARE} has a unique symmetric positive definite solution $P$,
  \item [(ii)] $\widecheck{U} \in \mathscr{F}$, where $\widecheck{U}(t,s)=u^*(t-s)$, $t\ge s \ge 0$, is an infinite horizon minimax observer.
  \item [(iii)] the minimax error is \(  \sigma(\widecheck{U},\ell) = \ell^{\top}F\LMAP^{\top}_g P \LMAP_g F^{\top}\ell \),
  \item [(iv)] $\widecheck{U}$ can be represented by the LTI system $(A_o,B_o,C_o,0)$, i.e.
     $\mobs{t}=C_o r(t)$ where $r\in AC([0,+\infty),\mathbb{R}^l)$ solves
   \begin{equation}
   \label{minimax:obs1}
     \begin{split}
     & \dot r(t) = A_or(t) + B_oy(t)  \mbox{\ and \ } r(0)=0,  \quad \mbox{ a.e. } \\
     \end{split}
   \end{equation}
 \end{itemize}
 \end{Theorem}
\begin{Remark}
   We stress that one does not need to solve~\eqref{minimax:obs1} several times in order to compute the minimax observer $\mobs{t}$ for various $\ell_i$, $i=1\dots s$.
Indeed, $A_o$ and $B_o$ do not depend on $\ell$,  since they depend only on the matrices $A_g,B_g,C_g,D_g$ of the associated stabilizable LTI $\mathscr{S}_g$ and on the
weight matrices $Q$ and $R$. In turn, $\mathscr{S}_g$ depends only on the matrices $F,A,H$, see Algorithm \ref{alg3}.
Hence, it is sufficient to define $C_o=L F\LMAP^T_g$ where $L=
   \left[\begin{smallmatrix}
     \ell_1\\\dots\\\ell_s
   \end{smallmatrix}\right]
$, provided $(F,A,H)$ is $\ell_i$-detectable, $i=1\dots s$.
 \end{Remark}
 \begin{Lemma}
 \label{obs:error}
  Assume $(F,A,H)$ is $\ell$-detectable and let $\widecheck{U}$ be the infinite horizon minimax observer from Theorem \ref{main:obs}. For any $(x,f,y,\eta) \in \mathscr{EE}([0,+\infty))$, the estimation error $e(t)=\ell^{\top}Fx(t)-\mobs{t}$ satisfies
  \begin{equation}
   \label{minimax:error}
    \begin{split}
    &  \dfrac{d\widetilde{r}}{dt} = (A_g-B_gK)^{\top}\widetilde{r}(t) + (C_g-D_gK)^{\top}\left[\begin{smallmatrix} f(t) \\ -\eta(t) \end{smallmatrix}\right]\,,  \\
    & \widetilde{r}(0)=(C_g-D_gK)^{\top} \left[\begin{smallmatrix} Fx(0) \\ 0 \end{smallmatrix}\right] \,, \\
    & e(t):=\ell^{\top}Fx(t)-\mobs{t}=\ell^{\top}F\LMAP^{\top}_g \widetilde{r}(t)\,.
    \end{split}
   \end{equation}
 In particular, if $f=0,\eta=0$, then $\lim_{t \rightarrow \infty} e(t)=0$.
 If $f,\eta$ are bounded, i.e. $\sup_{t > 0} \{f^\top(t)f(t),\eta^\top(t)\eta(t)\} <+\infty$, then the estimation error $e(t)$ is bounded,
 i.e. $\sup_{t > 0} e^\top(t)e(t) <+\infty$. If $f \in L^2([0,+\infty), \mathbb{R}^n)$, $\nu \in L^2([0,+\infty), \mathbb{R}^p)$, then $e \in L^2([0,+\infty),\mathbb{R})$.
 \end{Lemma}
 This lemma reveals the following important points: (i) the infinite horizon minimax observer from Theorem \ref{main:obs} behaves like an asymptotic observer in the absence of noises; (ii) if the noise signals are of bounded energy (i.e. they are in $L^2([0,+\infty)$), then the estimation error is of bounded energy; (iii) if the noise signals are just bounded, then the estimation error is bounded too, i.e., the upper bound of the error is proportional to the upper bound of the noises.
The discussion above allows us to formulate an algorithm for computing an infinite horizon state estimate (see Algorithm \ref{alg4}).  
The Matlab function
$DAEInfHorionsObs$ in the file $main\_TA3.m$ 
in the supplementary material of this report implements Steps \ref{alg4:step1} -- \ref{alg4:step2} Algorithm \ref{alg4}. 
\begin{algorithm}[H]
\caption{Infinite horizon minimax state estimate $\mobs{}$ \label{alg4}}
\renewcommand{\algorithmicrequire}{\textbf{Input:}}
\renewcommand{\algorithmicensure}{\textbf{Output:}}
\begin{algorithmic}[1]
\REQUIRE $(F,A,H,\ell)$ and the observed output $y \in L_{loc}^2([0,+\infty), \mathbb{R}^p)$
\ENSURE Infinite horizon minimax state estimate $\mobs{}$.
\STATE 
\label{alg4:step1}
Use Algorithm \ref{alg3} to compute the  stabilizable LTI system $\mathscr{S}_g=(A_g,B_g,C_g,D_g)$ associated with \eqref{dae:dual} and the state map $\LMAP_g$.
\STATE Compute the solution $P$ of \eqref{ARE} and compute the feedback matrix $K$ from \eqref{ARE}. In Matlab, the command $care(A_g,B_g, C_gSC_g^T, C_gSD_g^T, D_g^TSD_g)$ can be used to calculate $P$ and $K$.
\STATE
\label{alg4:step2}
 Set
\(
     A_o=(A_g-B_gK)^{\top},  ~B_o=(\hat{C}_u-\hat{D}_uK)^{\top}, ~ C_o=\ell^{\top}F\LMAP^T_g\,.
\)
\STATE
Solve the differential equation~\eqref{minimax:obs1} to find $r(t)$.
\STATE
  Return $\mobs{}=C_o r$.
\end{algorithmic}
\end{algorithm}

\section{Numerical example}
\label{sec:numerical-example}
Consider the following infinite-dimensional system with outputs, determined by the partial differential equation (PDE) of heat transfer:
 \begin{equation}
 \label{heat:eq1}
 \begin{split}
  & \dfrac{dV(t)}{dt} = \mathcal{A}V(t) +u(t)\,, \quad  y=\mathcal{C}V+w(t)\,,
 \end{split}
 \end{equation}
 where $V(t)$ and $u(t)$ take values in the space $H:=L^2([-1,1],\mathbb{R})$, i.e., $V(t)$ is itself a function from $[-1,1]$ to $\mathbb{R}$, and $\mathcal{A}v = -c\dfrac{d^2}{d\upsilon^2} v$ for all $v \in H$ such that $\mathcal{A}v\in H$. Moreover, $w \in L^1_{loc}([0,+\infty),\mathbb{R}^p)$, $\mathcal{C}:H \rightarrow \mathbb{R}^p$
is given by: $\mathcal{C}\phi=\int_{-1}^{1} g(\upsilon)\phi(\upsilon) d\upsilon$ for some $g \in L^2([-1,1],\mathbb{R}^p)$, $\phi \in H$.
We say that a pair $(V,u,y,w)$ solves~\eqref{heat:eq1} if
$V:[0,+\infty) \rightarrow H$ is a Frechet differentiable function such that $\mathcal{A}V(t)$ is defined, and \eqref{heat:eq1} holds a.e., and $V(t)(-1)=V(t)(1)=0$.
Physically, $V(t)(z)$ represents the temperature at time $t \in [0,+\infty)$ at point $z \in [-1,1]$. The constant $c$ in the definition of the operator $\mathcal{A}$
is the thermal diffusivity, in this example it is taken to be $c=0.033$ $\frac{\mathrm{m}}{\mathrm{s}^2}$.
From the control perspective, $V$ is a state, $u,w$ are noises and $y$ is a measured output.
We take $g(\upsilon)=(\sin(\pi n_1 \upsilon), \sin(\pi n_2 \upsilon)))^T$ and $p=2$, for other choices of $g$ and $p$ the discussion is similar.
Based on $y$ and assuming that $u$ and $w$ satisfy
$\int_0^{\infty} \|w(t)\|^2_{R^2}+ \|u(t)\|^2_H dt \le 1$, $\|.\|_H$ denotes the standard norm on $H$,
we would like to estimate $V(t)$ from $y$.

We recast the problem of estimating $V$ as the minimax estimation problem for DAEs as follows (see~\cite{ZhukSISC13} for the details): consider the basis $\{\phi_k\}_{k=0}^{\infty}$ of $H$, where $\phi_k=P_{k+1}-P_{k}$, $k \ge 1$, with $P_i$ denoting the $i$th Legendre polynomial, $i \ge 0$. Fix an integer $N > 0$ and define
\begin{equation*}
  \begin{split}
  & \forall \psi \in H: P_N(\psi)=\begin{bmatrix} \langle\psi,\phi_1\rangle, & \ldots, & \langle\psi,\phi_N\rangle \end{bmatrix}^T \\
  & \forall z=(z_1,\ldots,z_N)^T \in \mathbb{R}^N: P_N^{+}(z)=\sum_{k=1}^{N} z_k \phi_k \\
   & \hat{M}_N=(\langle\phi_i,\phi_j\rangle)_{i,j=1}^{N}, ~ A_N=(\langle\mathcal{A} \phi_i,\phi_j\rangle)_{i,j=1}^{N} \\
   & M_N = \Lambda_N (\langle\phi_i,\phi_j\rangle)_{i,j=1}^{N}, ~ \Lambda_N=(\frac{2i+1}{2} \delta_{i,j})_{i,j=1}^{N} \\
   & C_N = \begin{bmatrix} P_N (\sin(\pi n_1 \upsilon)), & P_N(\sin(\pi n_2 \upsilon)) \end{bmatrix}^T ,  \\
  \end{split}
 \end{equation*}
 where $\langle\cdot,\cdot\rangle$ is the standard inner product in $H$ and $\delta_{i,j}$ is the Kronecker delta symbol for all $i,j\ge 1$.
Consider the DAE \eqref{eq:dae} with  
 \[
  \begin{split}
   & F=\begin{bmatrix} M_N,  & 0  \\ 0 & 0 \\ 0 & 0 \end{bmatrix}, \quad  A=\begin{bmatrix} \Lambda_N A_N,  &  \Lambda_N \\ I_N & 0 \\ 0 & \begin{bmatrix} I_{N_u} & 0 \end{bmatrix} \end{bmatrix}, \quad H=\begin{bmatrix} C_N^T \\ 0  \end{bmatrix}^T
   \end{split}.
\]
Assume $(V,u,y,w)$ solve ~\eqref{heat:eq1}. 
From \cite{ZhukSISC13} it follows that the tuple $(x,f,y,\nu)$ is a solution of \eqref{eq:dae}, provided $a=\hat{M}_{N}^{-1}P_N V$,  $e_m=P_N\mathcal{A}\mathbf{e}$, $\mathbf{e}=V-P_N^{+}P_NV$, and $x=(a^T,e_m^T)^T$, and $\nu=\mathcal{C}\mathbf{e}+w$, where $f_1=\Lambda_N P_N u$, $f=\begin{bmatrix} f_1^T, & f_2^T, & f_3^T \end{bmatrix}^T$ for some $f_2 \in L^2_{loc}([0,+\infty),\mathbb{R}^N)$, $f_3 \in L^2_{loc}([0,+\infty),\mathbb{R}^{N_u})$. Moreover, if $V$ and its partial derivatives are bounded (which is the case for physically meaningful solutions), then by \cite{ZhukSISC13}, there exist $Q > 0$, $R > 0$ such that
$\int_0^{\infty} f^T(t)Qf(t)+\nu^T(t)R\nu(t)dt \le 1$. The intuition behind this is as follows:
the component $a(t)$ of $x(t)$ represents the first $N$ coordinates of $V(t)$ w.r.t the basis $\{\phi_k\}_{k=1}^{\infty}$,
and $e_m$ models the impact of $\mathbf{e}$, the error of approximating $V$ by its projection onto the first $N$ basis functions $\{\phi_k\}_{k=1}^N$, onto $a(t)$. For this reason, $e_m$ is included into the state of DAE. The latter allows to incorporate the projection error $\mathbf{e}$ into the estimation process. The smaller $e_m$ is, the closer the behavior of \eqref{eq:dae} is to that of \eqref{heat:eq1}.

The algebraic equations $a+f_2=0$ and $\left[\begin{smallmatrix} I_{N_u} & 0 \end{smallmatrix}\right]e_m+f_3=0$ together with $\int_0^{\infty} f^TQf(t) + \nu^T(t)R\nu(t)dt$ imply that
the $L^2$ norm of $a$ and of the first $N_u$ components of $e_m$ is uniformly bounded, i.e. these functions live in an ellipsoid.
Note that the expressions for $Q,R$ provided by \cite[Eq. (3.12) and (3.13)]{ZhukSISC13}  are very conservative, and, in practice, $Q,R$ are found by tuning.

The minimax observer for \eqref{eq:dae} can be used to estimate the state of \eqref{heat:eq1} as follows. Assume that \eqref{eq:dae} is $\ell$-detectable for $\ell \in \{\ell_i\}_{i=1}^{N_u}$, where $e_i$ is the $i$th standard basis vector in $\mathbb{R}^{2N+N_u}$. Let $\widecheck{U}_i$ denote the minimax observer corresponding to $\ell=\ell_i$ according to Theorem \ref{main:obs}, and let
$\mathscr{\vec{O}}_{N_u}(y)=(\mathscr{O}_{\widecheck{U}_1}(y),\ldots, \mathscr{O}_{\widecheck{U}_{N_u}}(y))^T$. Define  $\hat{\vec V}_{N_u}(y)(t)=P_{N_u}^{+}M_{N_u}^{-1}\mathscr{\vec O}_{N_u}(y)(t)$.
If $(V,u,y,w)$ is a solution of \eqref{heat:eq1}, then, as noted above, the tuple $(x,f,y,\nu)$ solves~\eqref{eq:dae}, and $P_{N_u}^{+}(\ell_1^TFx(t),\ldots,\ell_{N_u}^TFx(t))=P_{N_u}^{+}P_{N_u}V(t)$, the projection of $V$ to the first $N_u$ basis functions. By noting that $\mathscr{O}_{\widecheck{U}_i}(y)$ is the minimax estimate of $\ell_i^TFx$, we derive that $\hat{\vec V}_{N_u}(y)(t)$ is the minimax estimate of the vector $P_{N_u}^{+}P_{N_u}V(t)$. This latter estimate is ``good'' if the minimax error $\sigma(\widecheck{U}_i,\ell_i)$ is ``small''.

For the simulations we chose $N=40$, $n_1=30$, $n_2=31$, $N_u=10$. Note that the entries of $\hat{M}_N$, and $A_N$ can be computed analytically (see~\cite[Example 7.2, p. 121]{SpectralMethod}) so that $A,F,H$ can be computed numerically.
The numerical values of the matrices $F,A,H$ can be found in the file $DAE\_matrices.mat$ (Matlab format) or
in the files $F\_matrix.txt$, $A\_matrix.txt$, $H\_matrix.txt$ (csv text formal) of the supplementary material
of this report.
If we check $\rank \begin{bmatrix} F^T & A^T & H^T \end{bmatrix}^T$, then
it follows that its rank is less than its number of columns. Hence, by Remark \ref{rem_impulse_obs}, there exists $\ell$ such that $(F,A,H)$ is not $\ell$-impulse observable for this $\ell$. Thus, the observer proposed in~\cite{Darouach2009} is inapplicable as it requires to solve a linear equation \cite[eq. (5c)]{Darouach2009} which has no solution if there exists an $\ell$ for which $(F,A,H)$ is not $\ell$-impulse observable . 
 Note that \cite{Xu2007} is not applicable either: by \cite[Remark 1]{Xu2007}, existence of a filter in the sense of \cite{Xu2007} implies that the DAE at hand is regular, which is not the case.
However, using Algorithm \ref{alg5} - \ref{alg6}, it can be checked numerically that $(F,A,H)$ is $\ell$-impulse observable and $\ell$-detectable (Def.\ref{d:l-obs}-\ref{d:l-detect}) only for $\ell=\ell_i=e_i$, $i=1,\ldots,N_u$, i.e we can estimate (with finite worst-case error!) only the first $N_u$ components of $Fx$.

To generate the outputs we took an input $u_{true}(t)(\upsilon)= u_1(t) \sin(\pi n_1 \upsilon)+u_2(t)\sin(\pi n_2 \upsilon)$ with $u_1(t)=10cos(5t)$ and $u_2(t)=10sin(3t)$, initial condition $V_{true}(0)(\upsilon)=0.2\sin(\pi n_1 \upsilon)+0.2\sin(\pi n_2\upsilon)$, and noise
$w_{true}(t)=0.01e^{-0.001\pi^2c t}cos(100t)$. Then we computed the exact solution
 $(V_{true}, u_{true}, y_{true},w_{true})$ of~\eqref{heat:eq1} by using eigen-basis of the operator $\mathcal{A}$, i.e., $V_{true}(t,\upsilon)=z_1(t)\sin(\pi n_1 \upsilon)+z_2(t)\sin(\pi n_2 \upsilon)$ and
$y_{true}(t)=(z_1(t),z_2(t))^T+w_{true}(t)$, $\dot z_i = -c n_1^2 \pi^2 z_i+u_i(t)$, $i=1,2$. As noted, $(V_{true}, u_{true}, y_{true},w_{true})$ yields a solution of~\eqref{eq:dae}, namely
$(x_{true},f_{true},y_{true},\mu_{true})$ and $Fx_{true}(t)=\sum_{i=1}^2 z_i(t)\Lambda_NP_N(sin(\pi n_i \upsilon))$.

 We applied Algorithm \ref{alg6} to check that the DAE at hand is $\ell$-detectable for all standard unit
vectors of $\mathbb{R}^{N_u}$. We used the Matlab function $IsDetectable$ from the file $main\_TA3.m$ in the supplementary material of this report to check $\ell$-detectability.
We then applied Algorithm~\ref{alg4} to compute infinite horizon minimax observer $\mathscr{\vec{O}}_{N_u}(y_{true})(t)$ for \eqref{eq:dae} to estimate $Fx_{true}(t)$ on $[0,5]$. To this end we set $Q=\mathrm{diag}(10^{-5}I_{N_u},10^{-3}I_{N-N_u}, 10^{-3}I_N, 0.1I_{N_u})$, $R=400I_{p}$. We used the Matlab function
$DAEInfiniteHorizonObs$ from the file $main\_TA3.m$ of the supplementary material to calculate the 
matrices $A_o,B_o,C_o$ of the  minimax observer.
These matrices can be found in the files $Ao\_matrix.txt, Bo\_matrix.txt,Co\_matrix.txt$ (csv text format) and
$Observer.mat$ (Matlab data file format) included in the supplementary material of this report.
The results of the intermediate steps of the computation can also be found in the supplementary material.
More precisely
the matrices of the LTI $\mathscr{S}=(A_a,B_a,C_a,D_a)$ associated with the dual system \eqref{dae:dual} 
and the state map $\LMAP$
can be found  in the csv text files $Aa\_matrix.txt$, $Ba\_matrix.txt$, $Ca\_matrix.txt$, $Da\_matrix.txt$, $LMAP\_matrix.txt$
of the supplementary material of this report. 
The matrices of the stabilizable associated LTI 
$\mathscr{S}_g=(A_g,B_g,C_g,D_g)$ and the state map $\LMAP_g$ can be found in the csv text files
$Ag\_matrix.txt$, $Bg\_matrix.txt$, $Cg\_matrix.txt$, $Dg\_matrix.txt$, $LMAPg\_matrix.txt$
of the supplementary material of this report. The solution $P$ of \eqref{ARE} and the gain $K$ can be found in
the csv text files $P\_matrix.txt$ and $K\_matrix.txt$ of the supplementary material.

\begin{figure}[t]
\includegraphics[scale=0.4]{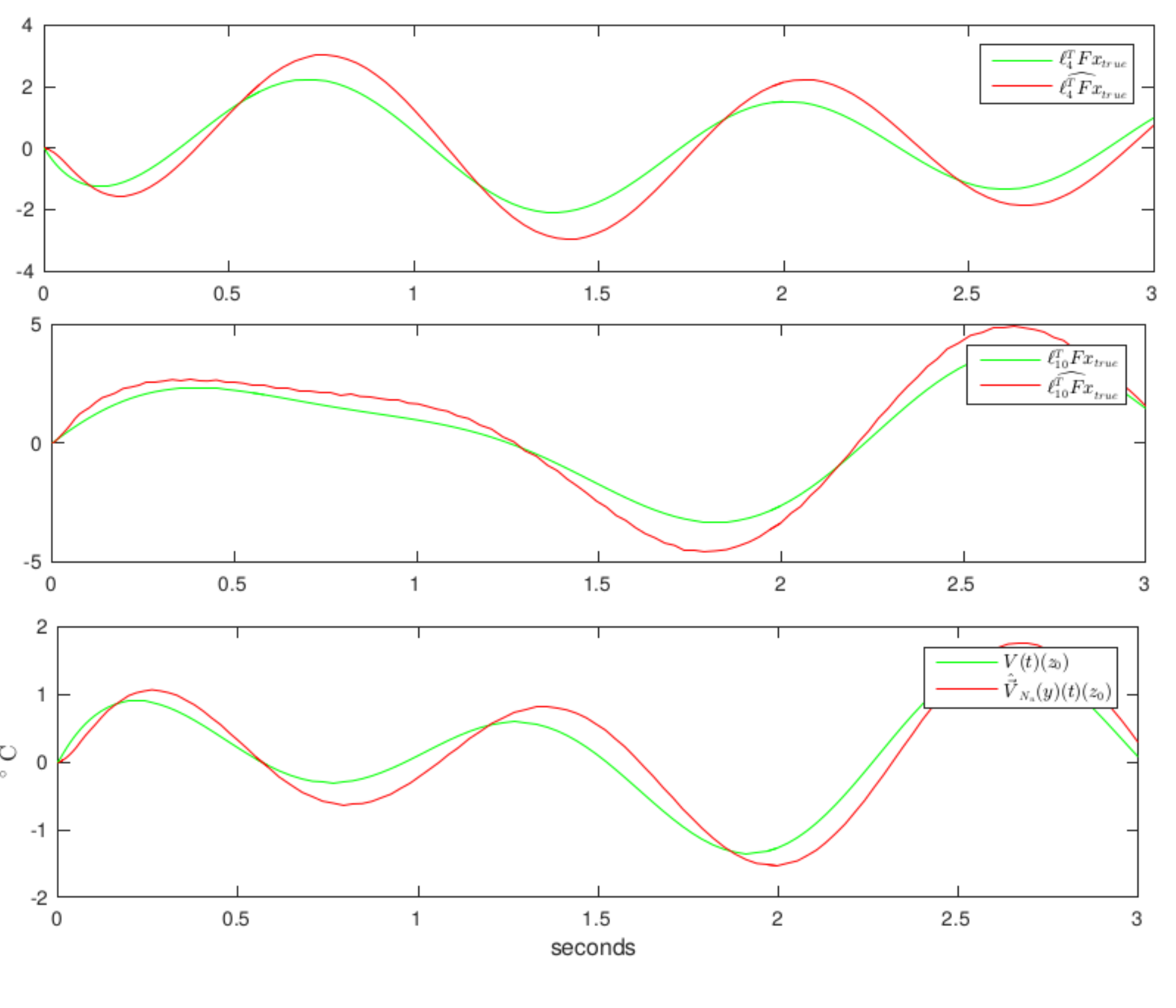}
\caption{Comparing $\widehat{\ell_i^TFx(t)}:=\mathscr{O}_{\widecheck{U}_{i}}(y_{true})(t)$ versus $\ell_{i}^TFx_{true}(t)$, $i=4,10$ and $V_{true}(t)(z_0)$  versus $\hat{\vec V}_{N_u}(y)(t)(z_0)$ for $z_0=-0.25$.
\label{fig}}
\end{figure}

It turns out that $\mathscr{\vec{O}}_{N_u}(y_{true})(t)$ is quite close to $Fx_{true}(t)$ component-wise.
 Specifically, Figure \ref{fig} shows that the minimax estimate $\mathscr{O}_{\widecheck{U}_i}(y_{true})(t)$ tracks the components $\ell^{T}_iFx_{true}(t)$, $i=4,10$ of DAE's state vector quite well.
Note that $Fx_{true}(t)$ represents the coordinates of the true solution $V_{true}$ with respect to the basis
$\{\phi_k\}_{k=0}^{\infty}$. Hence, it does not itself represent a physical quantity and therefore it has
no measurement unit.

We also compared $\hat{\vec V}_{N_u}(y)(t)$, the minimax estimate of the vector $P_{N_u}^{+}P_{N_u}V_{true}(t)$ against $V_{true}(t)$, and they were in a good agreement. In particular, on Figure \ref{fig} we present the plots of $t \mapsto V(t)(z_0)$ versus $t \mapsto \hat{\vec V}_{N_u}(y)(t)(z_0)$ for $z_0=-0.25$. Note that
  $\hat{\vec V}_{N_u}(y)(t)(z_0)$ represents an estimate of a physical quantity, namely, the temperature
  at time $t$ at point $z_0$.

\section{Proofs}
\label{sect:proof}
\begin{IEEEproof}[Proof of Theorem \ref{p:1:theo}]
Let us prove that
the set $\{(q,u) \in \mathscr{D}_{\ell}(I) \mid u=\delta_{t_1}(U)\}$ is not empty, if and only if $\sigma(U,\ell,t_1)<+\infty$.
This will also prove $\mathbf{(i)}\Leftrightarrow\mathbf{(ii)}$.
Take $w\in L^2(I)$ such that $F^{\top}w$ is absolutely continuous and $F^{\top}w(t_1)=F^{\top}\ell$. Define $
b(w,U)(t):=\dfrac{dF^{\top}w(t)}{dt}+A^{\top} w(t)-H^{\top}U(t)$. Clearly, $b(w,U) \in L^2(I)$. Now, let $(x,f,y,\eta)\in \mathscr{EE}(I)$. Then, integrating by parts (see~\cite[F.(2.1)]{Zhuk2012AMO}), we find:
\begin{equation}
\label{eq:lxobs}
\begin{split}
& \ell^{\top}Fx(t_1)-\obs{t_1} = (Fx(0))^{\top}(F^+)^{\top}F^{\top}w(0) + \\
 & + \int_0^{t_1} (w^{\top}(t)f(t)-U^{\top}(t)\eta(t) +b^{\top}(w,U)(t)x(t)) dt\,.
\end{split}
\end{equation}
Let $L$ denote a differential operator associated to~\eqref{eq:dae}: \((Lx)(t)=\bigl(Fx(0),\frac{d(Fx)}{dt}-Ax(t)\bigr)\) with domain $\mathscr D(L):=\{x\in L^2(I) \mid Fx \in AC(I), \frac{dFx}{dt} \in L^2(I) \}$.
Define the adjoint of $L$, namely $L'(g_0,g)(t) = -\frac{dF^{\top}g}{dt} - A^{\top}g\,, (g_0,g)\in\mathscr D(L')$,
 where the domain  $\mathscr D(L')$  of $L'$ is given by: \[
\begin{split}
& \mathscr D(L'):=\{(g_0,g) \in \mathbb{R}^{m} \times L^2(I)  \mid  F^{\top}g \in AC(I), \frac{dF^{\top}g}{dt} \in L^2(I)\,, \\
& F^{\top}g(t_1) = 0,  g_0 = {F^+}^{\top}F^{\top}g(0)+d,F^\top d=0\}\,.
\end{split}
\]
Now we claim that $\sigma(U,\ell,t_1)<+\infty\Leftrightarrow L'(g_0,g) = b(w,U)$ for some $(g_0,g)\in \mathscr D(L')$. Note that equality $L'(g_0,g) = b(w,U)$ implies that, by definition of $L'$ and $w$, $q:=w+g$ verifies $\frac{d(F^{\top}q)}{dt}=-A^{\top}q+H^{\top}U$ and $F^{\top}q(t_1)=F^{\top}\ell$, so that $(\delta_{t_1}(q),\delta_{t_1}(U)) \in \mathscr{D}_{\ell}(I)$, and vice versa, if $(\delta_{t_1}(q),\delta_{t_1}(U)) \in \mathscr{D}_{\ell}(I)$ then $g:=q-w$,
$g_0= {F^+}^{\top}F^{\top}g(0)$ verifies $L'(g_0,g) = b(w,U)$. Thus, it suffices to prove the above claim. To this end, recall the definition of $\mathscr{ E}(t_1)$ (Remark~\ref{r:rho_norm_Fx}). We note that $\sup_{(Fx(0),f,\eta)\in\mathscr E(t_1)} (Fx(0))^{\top}(F^+)^{\top}F^{\top}w(0) + \int_0^{t_1} w^{\top}(t)f(t)dt<+\infty$ by Cauchy-Swartz inequality in $L^2(I)$. Thus, by~\eqref{eq:error:fin} and~\eqref{eq:lxobs}, we have that $\sigma(U,\ell,t_1)<+\infty\Leftrightarrow s(b,U)<+\infty$, where \[
s(b,U):=\sup_{\{(x,\eta):(Lx,\eta)\in\mathscr E(t_1)\}}\int_0^{t_1}(b^{\top}(w,U)(t)x(t)-U^{\top}(t)\eta(t)) dt\,.
\]
To compute $s(b,U)$ we recall from convex analysis~\cite{Rockafellar1970} that the support function $s_E(x):=\sup_{v\in E}x^\top v $ of the ellipsoid $E:=\{v:v^\top Q_0 v\le 1\}$ equals to $(x^\top Q_0^{-1} x)^\frac 12$. Analogously, the support functional $s_{\mathscr E(t_1)}(g_0,g,v):=\sup_{(x_0,f,\eta)\in\mathscr E(t_1)}g_0^\top x_0+\int_0^{t_1}(g^\top f + v^\top\eta) dt$ of the ellipsoid $\mathscr E(t_1)$ equals to $(\rho(g_0,g,v,t_1,Q_0^{-1},Q^{-1},R^{-1}))^\frac 12$. Now we note that $s(b,U)$ equals to the support function of the pre-image of $\mathscr E(t_1)$ with respect to the operator $(x,\eta)\mapsto(Lx,\eta)$. The latter can be computed by using Young-Fenhel conjugate (see~\cite[F.(2.5),Lemma 2.2]{Zhuk2012AMO}):
\[
\begin{split}
&s^2(b,U) = \inf_{g_0,g,v}\{\rho(g_0,g,v,t_1,Q_0^{-1},Q^{-1},R^{-1}) \mid
(g_0,g) \in \mathscr{D}(L'), \\
& v \in L^2(I),
\bigl[\begin{smallmatrix}
L'(g_0,g) \\
v
\end{smallmatrix}\bigr]=\bigl[\begin{smallmatrix}
b(w,U) \\
U
\end{smallmatrix}\bigr]\}\,.
\end{split}\]
This equality implies that $s(b,U)<+\infty\Leftrightarrow L'(g_0,g) = b(w,U)$ for some $(g_0,g)\in \mathscr D(L')$. 

Let us now prove ~\eqref{eq:sigmaUt1}. If $\sigma(U,\ell,t_1)=+\infty$, then from the discussion above
both sides of \eqref{eq:sigmaUt1} are $+\infty$.
Assume $\sigma(U,\ell,t_1) <  +\infty$.
Take $(p,u)\in \mathscr{D}_\ell(I)$ such that $U=\delta_{t_1}(u)$. Set $w:=\delta_{t_1}(p)$. Then $b(w,U)=0$ in~\eqref{eq:lxobs} and so we get:
\begin{equation}
  \label{eq:sigmaUc}
 \begin{split}
 & (\sigma(U,\ell,t_1))^\frac 12 = s_{\mathscr E_1(t_1)}(w_0,w,U):=\\
 & \sup_{(x_0,f,\eta) \in \mathscr E_1(I)} \{ w_0^{\top}x_0 +
   \int_0^{t_1} (w^{\top}(t)f(t)-U^{\top}(t)\eta(t)) dt\}\,,
 \end{split}
\end{equation}
where $w_0:=(F^+)^{\top}F^{\top}w(0)$, $\mathscr E_1(t_1)$ is the set of all $(x_0,f,\eta)\in \mathscr E(t_1)$ such that $Lx = (x_0,f)$ for some $x\in\mathscr D(L)$,
and $s_{\mathscr E_1(t_1)}$ is the support function of $\mathscr E_1(t_1)$.
Thus, to compute $\sigma(U,\ell,t_1)$ it suffices to compute
right hand side of~\eqref{eq:sigmaUc}.
 From ~\cite[F.(2.6),Lemma 2.2]{Zhuk2012AMO}:
\begin{equation}
  \label{eq:czU}
 \begin{split}
& s_{\mathscr E_1(t_1)}^2(w_0,w,U)= \\
& \inf_{(g_0,g)\in\ker L'}\rho(w_0+g_0,w+g,U,t_1,Q_0^{-1},Q^{-1},R^{-1})\,,
\end{split}
\end{equation}
where $\ker{L'}=\{(g_0, g)\in\mathscr D(L'):L'(g_0,g)=0\}$.
 We claim that the right-hand side of~\eqref{eq:czU} equals to $\inf_{(q,v) \in \mathscr{D}_{\ell}(I), v=\delta_{t_1}(U)}J(q,v,t_1)$.  Indeed, since $g_0 = {F^+}^{\top}F^{\top}g(0)+P\tilde d$, $\tilde d\in\mathbb R^m$, \(P=(\Id_m-(F^{\top})^+F^{\top})\) it follows that we can fix $g$ and eliminate $\tilde d$ by minimizing the 1st term of $\rho$ w.r.t. $\tilde d$. We get:
\begin{equation*}
  \begin{split}
  & \inf_{\{g_0 \mid (g_0,g) \in \ker L'\}} \|Q_0^{-\frac{1}{2}}((F^+)^{\top}F^{\top}w(0)+g_0)\|^2= \\
  & =\inf_{\tilde d \in \mathbb{R}^m}  \|Q_0^{-\frac 12}((F^+)^{\top}F^{\top}(w(0)+g(0))-P\tilde d)\|^2 \\
  & =\|Q_0^{-\frac 12} ((F^+)^{\top}F^{\top}(w(0)+g(0))-P\hat d)\|^2= \\
  & = (w(0)+g(0))^{\top}\bar{Q}_0(w(0)+g(0))\,,
\end{split}
\end{equation*}
where $\hat d:=(PQ_0^{-1}P)^{+}PQ_0^{-1}(F^+)^{\top}F^{\top}(w(0)+g(0))$. Now, to prove the above claim we note that \(\{(q,v)\mid v=u,(q,v)\in \mathscr{D}_\ell(I)\}=\{(\delta_{t_1}(g+w),u)\mid (g_0,g)\in\ker{L'}\}\). In other words, if $(q,u)\in \mathscr{D}_\ell(I)$ then $q$ is a sum of $p=\delta_{t_1}(w)$ and $\delta_{t_1}(g)$ provided $L'(g_0,g)=0$. This allows us to write:
\begin{equation*}
\begin{split}
  & \inf_{\{g\mid(\hat g_0,g)\in\ker L'\}} [(w(0)+g(0))^{\top}\bar{Q}_0(w(0)+g(0))+ \\
   & \int_0^{t_1}(w(t)+g(t))^{\top}Q^{-1}(w(t)+g(t))dt]  =  \inf_{(q,u) \in \mathscr{D}_{\ell}(I), u=\delta_{t_1}(U)} [ \\
   &(q(t_1))^{\top}\bar{Q}_0q(t_1) +\int_0^{t_1}q^{\top}(t)Q^{-1}q(t)dt]\,,
\end{split}
\end{equation*}
where $\hat g_0 = {F^+}^{\top}F^{\top}g(0)+P\hat d$. Thus, the right-hand side of~\eqref{eq:czU} equals $\inf_{(q,v) \in \mathscr{D}_{\ell}(I), v=\delta_{t_1}(U)}J(q,v,t_1)$. The latter and~\eqref{eq:sigmaUc} prove~\eqref{eq:sigmaUt1}.
%
To prove $\mathbf{(ii)}\Leftrightarrow\mathbf{(iii)}$ we note that, by Definition~\ref{d:finhorobs}, $\mathbf{(iii)}$ implies that $\sigma(\widehat U,\ell,t_1)<+\infty$ and then $\mathbf{(ii)}$ follows. Now, assume that $\mathbf{(ii)}$ holds and $(q^*,u^*) \in \mathscr{D}_{\ell}(I)$ denotes a minimizer of $J$ over the affine set $\mathscr{D}_{\ell}(I)$.
Then, by~\eqref{eq:sigmaUt1}, $\widehat{U}=\delta_{t_1}(u^*)$ is a minimizer for $\sigma(U,\ell,t_1)$ and so $\mathbf{(iii)}$ holds by Definition~\ref{d:finhorobs}.
\end{IEEEproof}
\begin{IEEEproof}[Proof of Theorem \ref{p:5}]
To prove $\mathbf{(i)}\Rightarrow\mathbf{(ii)}$ we assume that $(F,A,H)$ is $\ell$-detectable. Then \eqref{dae:dual} is behaviorally stabilizable from $\ell$ (according to the terminology of~\cite{AutomaticaPaper}) and vice-versa. Then by~\cite[Theorem 6]{AutomaticaPaper} there exists
$(q^*,u^*) \in \mathscr{D}_{\ell}([0,+\infty))$ such that~\eqref{eq:dual_LQ_inf} holds. Let us prove that $\widecheck{U}(t,s)=u^*(t-s)$ verifies~\eqref{problem:obs:eq2}. To this end, from \eqref{eq:sigmaUt1} it follows that for any $u \in L^2([0,\tau],\mathbb{R}^p)$:
\begin{equation}
  \label{eq:minIminJ}
\begin{split}
& \sigma(u,\ell,\tau) 
\ge    \inf_{(q,v) \in \mathcal{D}_{\ell}([0,\tau]) } J(q,v,\tau)\,.
\end{split}
\end{equation}
Take any $U \in \mathscr{F}$. Then by using~\eqref{eq:minIminJ} and~\eqref{eq:dual_LQ_inf} and \eqref{problem:obs:eq2.1}
 we get:
\begin{equation}
  \label{eq:barU}
 \begin{split}
& \sigma(U,\ell)  
 \ge \limsup_{\tau \rightarrow \infty} \inf_{(q,v) \in \mathcal{D}_{\ell}([0,\tau]) } J(q,v,\tau) =  J^*_\infty\,.
\end{split}
\end{equation}
Apply \eqref{eq:sigmaUt1} to $\widecheck U(\tau,\cdot)$. Then \( \sigma(\widecheck{U}(\tau,\cdot),\ell,\tau) = \inf_{(q,v) \in \mathcal{D}_{\ell}([0,\tau]), v={u^*}|_{[0,\tau]} } J(q,v,\tau)\le J(q^*,u^*,\tau)\), so that \( \sigma(\widecheck{U},\ell)\le \limsup_{\tau \rightarrow \infty} J(q^*,u^*,\tau)=J^*_\infty\). This and~\eqref{eq:barU} implies that $\sigma(\widecheck{U},\ell)\le J^*_\infty \le \sigma(U,\ell)$. Since $J^*_\infty<+\infty$ by~\eqref{eq:dual_LQ_inf} it follows that $\sigma(\widecheck{U},\ell)<+\infty$ and so $\widecheck U$ verifies~\eqref{problem:obs:eq2}.

To prove $\mathbf{(ii)}\Rightarrow\mathbf{(i)}$ let $\widecheck{U} \in \mathscr{F}$ be an infinite horizon minimax observer. Then, by Definition~\ref{d:finhorobs},
 $\sigma(\widecheck{U},\ell):=\limsup_{t \rightarrow \infty}\sigma(\widecheck{U}(t,\cdot),\ell,t)<+\infty$. Hence, there exists $\tilde t>0 $ such that for all $t > \tilde t$: $\sigma(\widecheck{U}(t,\cdot),\ell,t)<+\infty$. Now, by \eqref{eq:sigmaUt1} of Theorem~\ref{p:1:theo} we get that   $\sigma(\widecheck{U}(t,\cdot),\ell,t)=\inf_{(q,u) \in \mathscr{D}_{\ell}([0,t]), u=\delta_{t}(\widecheck{U}(t,\cdot))} J(q,u,t)$ for all $t > \tilde t$, and so $\limsup_{t \rightarrow \infty} \inf_{(q,u) \in \mathscr{D}_{\ell}([0,t]), u=\delta_{t}(\widecheck{U}(t,\cdot))} J(q,u,t)<+\infty$. The latter and~\cite[Theorem 6]{AutomaticaPaper} implies that~\eqref{dae:dual} is behaviorally stabilizable from $\ell$ that, as noted above, implies $\mathbf{(i)}$.
\end{IEEEproof}
\begin{IEEEproof}[Proof of Theorem \ref{t:observer-fin-hor}]
Assume that $v,u^*$ are defined by~\eqref{dual:finhor:sol}. Set $q^*(s)=(C_s+D_sK(t_1-s))v(s)$, $s \in I$. From \cite[Theorem 5]{AutomaticaPaper}
it then follows that $J(q^*,u^*,t_1)=\inf_{(q,u) \in \mathscr{D}_{\ell}(I)} J(q,u,t_1)$ and
$J(q^*,u^*,\tau_1)=\ell^{\top}F\LMAP^{\top} P(t_1)\LMAP F^{\top}\ell$.
From Theorem \ref{p:1:theo} it follows that $\widehat{U}$ is indeed the minimax observer on $I$ and \( \sigma(\widehat{U},\ell,t_1)=J(q^*,u^*,t_1)=\ell^{\top}F\LMAP^{\top} P(t_1)\LMAP F^{\top}\ell. \) Finally, let $r$ be the solution of \eqref{finhor:eq1}. Then
$\frac{d}{dt} (r^{\top}(t)v(t_1-t))=y^{\top}(t)u^*(t_1-t)=y^{\top}(t)\widehat{U}(t)$, and by integrating both sides,
$\mathcal{O}_{\widehat{U}}(t_1,y) =\ell^{\top}F\LMAP^{\top}r(t_1)$.
\end{IEEEproof}
\begin{IEEEproof}[Proof Theorem \ref{main:obs}]
 If $(F,A,H)$ is $\ell$-detectable, then \eqref{dae:dual} is behaviorally stabilizable from $\ell$ (according to the terminology of~\cite{AutomaticaPaper}) and vice-versa. Hence, the statement of $\mathbf{(i)}$ follows from~\cite[Lemma 2]{AutomaticaPaper}.
From \cite[Theorem 6]{AutomaticaPaper} it follows that if $v^* \in AC([0,+\infty),\mathbb{R}^l)$ is a solution of
 $\dot {v^*} = (A_g-B_gK)v^*\,, \quad v^*(0)=\LMAP_g(F^{\top}\ell)$ and $(q^*,u^*) \in L^2_{loc}([0,+\infty),\mathbb{R}^m) \times L^{2}_{loc}([0,+\infty),\mathbb{R}^p)$
 are such that $\left[
   \begin{smallmatrix}
     q^*\\ u^*
   \end{smallmatrix}\right]=(C_g-D_gK)v^*$, then $(q^*,u^*) \in \mathscr{D}_{\ell}([0,+\infty))$ and $J^*_\infty = \limsup_{t \rightarrow \infty} J(q^*,u^*,t) = \limsup_{t \rightarrow \infty} \inf_{(q,u) \in \mathscr{D}_{\ell}([0,t])} J(q,u,t)$ and $J^*_\infty = \ell^{\top}F\LMAP^{\top}_g P \LMAP_g F^{\top}\ell$. From $\mathbf{(ii)}$, Theorem~\ref{p:5} it then follows that $\mathbf{(ii),(iii)}$ hold true. Now, to prove $\mathbf{(iv)}$ we note that
  from the above discussion $u^{*}(t)=(C_{u}-D_uK)e^{(A_g-B_gK) t}\LMAP_gF^{\top}\ell$ and thus $(A_o,B_o,C_o)$ represents $\widecheck{U}$.
\end{IEEEproof}
\begin{IEEEproof}[Proof of Lemma \ref{obs:error}]
   Let  $(q^*,u^*) \in \mathcal{D}_{\ell}([0,+\infty))$ be such that $({q^*}^\top,{u^*}^\top)^\top = (C_g-D_gK)p$, $\dot p = (A_g-B_gK)p$, $p(0)=\LMAP_g(F^{\top}\ell)$.
  Define $g:[0,t] \mapsto (Fx)^\top(s)(F^\top)^{+}F^\top q^*(t-s)=x^\top(s)F^\top q^*(t)$. Since $q^*$ and $Fx$ are both absolutely continuous, $g$ is absolutely continuous. Note that $Ax(s)+f(s)$ is in the range of $F$ a.e. and so $FF^{+}(Ax(s)+f(s))=Ax(s)+f(s)$ a.e. Hence,
 \[ \frac{dg(s)}{ds}= f^{\top}(s)q^*(t-s)+(y(s)-\eta(s))^{\top}u^*(t-s)\, a.e. \]
  Integrate both sides of the above equality and use $\widecheck{U}(t,s)={u^*}(t-s)$:
  \begin{equation}\label{obs:error:eq1}
   g(t)- \mobs{t} = g(0)+ \int_0^{t} \left[\begin{smallmatrix} {q^*}^{\top}(t-s) & {u^*(t-s)} \end{smallmatrix}\right]\left[\begin{smallmatrix} f(s) \\ -\eta(s) \end{smallmatrix}\right]ds\,.
  \end{equation}
  Note that $g(t)=\ell^{\top}Fx(t)$, and hence the left-hand side of \eqref{obs:error:eq1} is the estimation error $e(t)$. It is easy to see that
  the right-hand side of \eqref{obs:error:eq1} equals $\ell^{\top}F\LMAP^{\top}_g \widetilde{r}(t)$ where $\widetilde{r}$ solves~\eqref{minimax:error}.
\end{IEEEproof}

\section{Conclusions}
\label{sec:conclusions}
We have presented an original solution of the minimax observer design problem for linear DAEs. It is ``application friendly'' as (i) it relies upon well developed tools (e.g. Riccati equations), (ii) it imposes easy to check necessary and sufficient conditions ($\ell$-detectability) on $(F,A,H)$ and (iii) it is optimal for $L^2$ inputs and provides bounded estimation error for merely bounded inputs.

\bibliographystyle{IEEEtran}


\end{document}